\documentclass[screen]{aomart}
\makeatletter
\fancypagestyle{firstpage}{%
  \fancyhf{}%
  \if@aom@manuscript@mode
    \lhead{\begin{picture}(0,0)%
        \put(-21,-25){\usebox{\@aom@linecount}}%
      \end{picture}}
  \fi
  \chead{}%
   \cfoot{\footnotesize\thepage}}%
\makeatother

\usepackage{cancel}
\usepackage{amsthm}
\usepackage{amsmath}
\usepackage{yfonts}
\chardef\bslash=`\\ 

\newtheorem[{}\it]{thm}{Theorem}
\newtheorem{cor}[thm]{Corollary}
\newtheorem{lem}[thm]{Lemma}
\newtheorem{prop}[thm]{Proposition}

\newtheorem[{}\it]{thmO}{Theorem}
\newtheorem{corO}[thmO]{Corollary}
\newtheorem{lemO}[thmO]{Lemma}
\newtheorem{propO}[thmO]{Proposition}

\theoremstyle{definition}
\newtheorem{defn}{Definition}
\newtheorem{rem}{Remark}

\newtheorem*[{}\it]{notation}{Notation}

\newtheorem{defnO}{Definition}

\newcommand{\CWF}[2][c]{\begin{tabular}[#1]{@{}c@{}}#2\end{tabular}}

\addresses

\title[Universal geometries of linear $2^\mathrm{nd}$ order ordinary diff. eqs.]{Universal geometries underpinning linear second order ordinary differential equations}
\author{{\L}ukasz Rudnicki}
\address{Faculty of Mathematics, Physics, and Informatics \& International Centre for Theory of Quantum Technologies, 
 University of Gdansk,
80-308 Gdansk, Poland}
\email{lukasz.rudnicki@ug.edu.pl} 
\orcid{0000-0001-8563-6101}


\copyrightnote{}

\begin{document}
\begin{abstract}
A deep relationship \cite{Rudnicki} between real linear second order ordinary differential equations $u''\left(x\right)+h\left(x\right)u\left(x\right)=0$, with differentiable $h(x)$, and two dimensional hyperbolic geometry is generalized in a multitude of ways. First of all, still in the real case, I present an equivalent relationship in which the hyperbolic geometry is replaced by a two dimensional (anti-)de Sitter geometry. As happened in the former case, I show that this equation everywhere admits a pair of linearly independent solutions locally expressed in terms of an arbitrary non-vertical geodesic curve in the discussed geometry. Moreover, I show that every solution of a corresponding Ricatti equation $ \Theta'\left(x\right)+\Theta^2\left(x\right)+h(x)=0$ obtained through a substitution $u'\left(x\right)=\Theta\left(x\right)u\left(x\right)$, still with differentiable $h(x)$, itself is a geodesic curve in the two dimensional (anti-)de Sitter geometry.

Next, after promoting $h(x)$ to a holomorphic function $h(z)$, I express two linearly independent solutions of $u''\left(z\right)+h\left(z\right)u\left(z\right)=0$ in virtually the same way as for the real scenario and hyperbolic geometry. In this case, the curves used to build the solutions are geodesic in a two dimensional complex Riemannian geometry of a sphere. Analogous results for the complex Ricatti equation follow. As before, this geometric interpretation is independent of the function $h(z)$, while the holomorphic metric assumes the same functional form as the hyperbolic metric discovered in \cite{Rudnicki}, with real coordinates just replaced by complex coordinates.

Finally, I show that  complex linear second order ordinary  differential equation in question is in a relationship with
four dimensional pseudo Riemannian K{\"a}hler-Norden geometry, equivalent to its relationship with the  complex Riemannian geometry. The added value of working with real geometry turns out to be that certain two dimensional submanifold of the K{\"a}hler-Norden manifold render the hyperbolic and the (anti-)de Sitter scenario, both relevant for the real equation. Thus, this last discussed geometry serves as the base model for all considered geometric underpinnings of linear second order ordinary differential equations.
\end{abstract}

\maketitle
\tableofcontents

\section{Introduction}
We consider a linear second order ordinary differential equation
\begin{equation}\label{Eq.mainC}
u''\left(z\right)+h\left(z\right)u\left(z\right)=0,
\end{equation}
with a complex independent variable $z$ and an arbitrary complex-valued function $h(z)$. In a prequel manuscript \cite{Rudnicki} I have shown there exists a generic connection between (\ref{Eq.mainC}), for $z$ being a real variable $x$ and $h(x)$ being a real-valued function of class $C^1\left(\mathbb{R}\right)$, and two dimensional hyperbolic geometry given by an $h$-dependent Riemannian metric with sectional curvature equal to $-1$. 

In this contribution I aim to extend this universal geometric underpinning to a complex scenario. Even though complex numbers double the real dimension of the problem, thus, two dimensional real geometry is not enough in general, in \fullref{Sec.}{Sec2}, I present the following results. First of all, two dimensional (anti-)de Sitter geometry (maximally symmetric Lorentzian geometry of constant sectional curvature) provides an alternative complex-valued description of the real scenario, with a benefit that the geodesic equation in this geometry can as well directly be reduced to the Ricatti equation. I use the term \textit{(anti-)de Sitter}, because de Sitter geometry in two dimensions only differs from its anti-de Sitter counterpart by a minus sign in front of its Lorentzian metric. Next, as the major result I find that holomorphic Riemannian geometry \cite{LeBrun1,LeBrun2,ganchev} allows to translate all technical findings from \cite{Rudnicki} to the general case of Eq. (\ref{Eq.mainC}), provided that $h(z)$ is holomorphic. A topologically distinct feature is that in holomorphic Riemannian geometry  spaces with a constant nonzero holomorphic sectional curvature are locally holomorphical isometric to a complex sphere \cite{ganchev}. By analogy with the hyperbolic case, our complex holomorphic metric has its holomorphic sectional curvature equal to $-1$, what corresponds to a complex 2-sphere of (complex) radius equal to imaginary unit $i$. Going back to real differential geometry, we find that all differential equations (\ref{Eq.mainC}) with holomorphic $h(z)$ correspond to a pseudo-Riemannian K{\"a}hler-Norden geometry \cite{Borowiec} in dimension four, described by a metric which is Einstein
. This geometry possesses two distinct two dimensional submanifolds  which are equal to both descriptions of the real scenario: hyperbolic from \cite{Rudnicki} and (anti-)de Sitter  from \fullref{Sec.}{subsec2.1}, respectively.

Following a few concluding  remarks in \fullref{Sec.}{Sec3}, proofs of all results are collected in \fullref{Sec.}{Sec4}. Whenever possible, both in the main text and in the proofs, I refer to the aforementioned hyperbolic scenario. To facilitate this approach, in \fullref{Sec.}{subsechyp}, I collect all necessary prerequisites. To avoid confusion, I enumerate definitions, theorems, etc., quoted from \cite{Rudnicki} (sometimes with adjustments), with "H" in front of their numbers referring to "hyperbolic".

\subsection{Real scenario and hyperbolic geometry}\label{subsechyp}
In \cite{Rudnicki} I studied a linear second order ordinary differential equation
\begin{equation}\label{Eq.main}
u''\left(x\right)+h\left(x\right)u\left(x\right)=0,
\end{equation}
with the independent variable $x$ assumed to be real, and with $h(x)$ being an arbitrary real-valued function of class $C^1$ (sometimes $C^2$). I introduced a two dimensional Riemannian geometry associated with Eq. (\ref{Eq.main}), found it is hyperbolic, discussed geodesic curves in this geometry and proved a local relation between geodesic curves in explicit form (with affine parameter eliminated) and a general solution of (\ref{Eq.main}). Below, I shall only list tailored versions of results which are essential for the main body of this manuscript. Proofs, accompanying statements, thorough contextualization and examples can all be found in the original paper. A few proofs are also repeated in \fullref{Sec.}{Sec5}.
\begin{defnO}[Definition 3 from \cite{Rudnicki}: Linear-$2^{\mathrm{nd}}$-order-ODE upper half plane] \label{DefM}
Let $\mathcal{H}_0 =\left\{\left(x,\varPhi\right)\in\mathbb{R}\times\mathbb{R}_{>0}\,|\,\varPhi^2\neq h(x)\right\}$. Let
\begin{equation}\label{metric}
    g_h=\frac{\left(h(x)-\varPhi^2\right)^2dx^2+d\varPhi^2}{\varPhi^2},
\end{equation}
be a metric on $\mathcal{H}_0$. A pair $\mathbb{M}_h = \left(\mathcal{H}_0, g_h\right)$ is a Riemannian manifold.
\end{defnO}
\begin{lemO}[Lemma 1 from \cite{Rudnicki}]
Let $h(x)$ be of class $C^2$. The metric $g_h$ in the whole domain $\mathcal{H}_0$ has constant sectional curvature equal to $-1$. Therefore, it everywhere describes hyperbolic geometry. \label{lem-hyp}
\end{lemO}

With the dot we  denote derivatives with respect to an affine parameter "$s$". In some local coordinates $\varsigma^k$, the geodesic equation \cite{Besse} $\dot{\varsigma}^j\nabla_j\dot{\varsigma}^i=0$ becomes $\ddot{\varsigma}^i+\Gamma_{\:\:\:jk}^{i}\left(\varsigma\right)\dot{\varsigma}^j\dot{\varsigma}^k=0$, where $\nabla$ is the unique Levi-Civita connection associated with $g_h$ and we use Einstein summation convention. 

\begin{defnO}[Definition 5 from \cite{Rudnicki}: Geodesic curve on $\mathbb{M}_h$ in explicit form] \label{defGEO} Let a geodesic curve on $\mathbb{M}_h$, $\mathbb{R}\ni s \mapsto \left(x(s),\varPhi(s)\right)\in \mathcal{H}_0$, be a solution of the geodesic equation with some initial conditions set for $s=0$. If
\begin{equation}\label{inverdom}
    \displaystyle\mathop{\forall}_{s\in\left]s_{-},s_{+}\right[} \;\; \dot{x}\left(s\right)\neq 0,
\end{equation} holds for some $s_-<0<s_+$, then we can re-express the geodesic curve through a substitution $\varPhi(s)= \Phi\left(x(s)\right)$ and  call the function $\Phi\left(x\right)$ a geodesic curve on $\mathbb{M}_h$ represented in explicit form. It exists on $\left]x_-,x_+\right[$, where $x_-<x_0<x_+$ and $x_0=x(0)$. If $\dot{x}_0 > 0$, then  $x_-=x(s_-)$ and $x_+=x(s_+)$, while if $\dot{x}_0 < 0$, then  $x_-=x(s_+)$ and $x_+=x(s_-)$.
\end{defnO}
\begin{rem}\label{RemI1}
The geodesic curve is subject to generic initial conditions: $x(0)=x_0$, $\dot{x}(0)=\dot{x}_0$, $\varPhi(0)=\varPhi_0$ and $\dot{\varPhi}(0)=\dot{\varPhi}_0$. Since in this way we can assure that $\dot{x}_0\neq 0$, the above construction always applies, at least locally.
\end{rem}

\begin{propO}[Proposition 2 from \cite{Rudnicki}]
Let $h(x)$ be of class $C^1$. A geodesic curve on $\mathbb{M}_h$ in explicit form $\Phi(x)$ obeys \label{Prop1}
\begin{equation}\label{ELexp}
 \Phi''\left(x\right) 
 =\frac{3\Phi^{2}\left(x\right)+h\left(x\right)}{\Phi^{2}\left(x\right)-h\left(x\right)}\frac{\left[\Phi'\left(x\right)\right]^{2}}{\Phi\left(x\right)}-\frac{h'\left(x\right)\Phi'\left(x\right)}{\Phi^{2}\left(x\right)-h\left(x\right)}+\frac{\Phi^{4}\left(x\right)-h^{2}\left(x\right)}{\Phi\left(x\right)}.
\end{equation}
We observe that Eq. (\ref{ELexp}) can  be written as
\begin{equation}\label{extradop}
 \Phi''\left(x\right) =\Lambda\left(x,\Phi(x)\right)\left[\Phi'(x)\right]^{2}+\Gamma_{\:\:\:xx}^{x}\left(x,\Phi(x)\right)\Phi'(x)-\Gamma_{\:\:\:xx}^{\varPhi}\left(x,\Phi(x)\right),
\end{equation}
where $\Lambda\left(x,\Phi(x)\right)=2\Gamma_{\:\:\:\varPhi x}^{x}\left(x,\Phi(x)\right)-\Gamma_{\:\:\:\varPhi\varPhi}^{\varPhi}\left(x,\Phi(x)\right)$.
\end{propO}
\begin{notation}[Coordinates and geodesic curves in explicit form] With italic font I denote the coordinate $\varPhi$ and its $s$-dependent variant $\varPhi(s)$ (which is the coordinate of the geodesic curve), while functions of the variable $x$ representing geodesic curves in explicit form $\Phi(x)$ are denoted with upright $\Phi$. For example, $\Gamma_{\:\:\:\varPhi\varPhi}^{\varPhi}\left(x,\Phi(x)\right)$ represents the $\Gamma_{\:\:\:\varPhi\varPhi}^{\varPhi}$ Christoffel symbol (coordinates) evaluated at $x$ and $\Phi(x)$ (a point on a geodesic curve in explicit form).
\end{notation}
\begin{rem}\label{RemI2}
    In \fullref{Sec.}{Sec5} we prove \fullref{Proposition}{Prop1} by substituting $\varPhi(s)= \Phi\left(x(s)\right)$ into the geodesic equation. Then, we find  that the geodesic equation implies (\ref{ELexp}) as long as $\dot{x}(s) \neq 0$, i.e., the equation locally governing $\Phi(x)$ decouples from the equation for $x(s)$. Thus, we can first find a solution of (\ref{ELexp}), and letting $x=x(s)$ substitute it back into the geodesic equation. After finding $x(s)$ we impose initial conditions and obtain the full trajectory. In other words, the condition  $\dot{x}(s) \neq 0$ holding on some $\left]s_{-},s_{+}\right[$, can be verified \textit{a posteriori}. As a consequence, this condition is enough to consider Eq. \ref{ELexp}, without referring to the inverse function $s(x)$ and its existence, the latter giving $\Phi(x)=\varPhi\left(s(x)\right)$. While in the current example it is a tautology because the inverse function is guaranteed, this observation will be important in a multi-variable scenario.  
\end{rem}

\begin{thmO}[Theorem 3 from \cite{Rudnicki}]
Solutions of all linear second order ordinary differential equations (\ref{Eq.main}), with any $h(x)$ of class $C^1$, can be locally expressed through geodesic curves in hyperbolic two dimensional geometry.  \label{thm-main2}
    
More precisely, let  $h(x)$ be of class $C^1$, let $x_0\in\mathbb{R}$ and let 
$\Phi(x)$ be of class $C^2$ on some interval $\left]x_-,x_+\right[$, such that $x_-<x_0<x_+$ and $h\!\left(x\right)\neq\Phi^{2}\!\left(x\right)$ on this interval. If
\begin{subequations}\label{specsol}
\begin{equation}\label{specsol1}
        u_\mathtt{top}(x)=\exp\left[\int_{x_0}^{x}\!\!d\xi\,\Phi\left(\xi\right)\frac{\Phi'\left(\xi\right)-\sqrt{\left[h\left(\xi\right)-\Phi^{2}\left(\xi\right)\right]^{2}+\left[\Phi'\left(\xi\right)\right]^{2}}}{h\left(\xi\right)-\Phi^{2}\left(\xi\right)}\right],
    \end{equation}
    \begin{equation}\label{specsol2}
        u_\mathtt{bot}(x)=\exp\left[\int_{x_0}^{x}\!\!d\xi\,\Phi\left(\xi\right)\frac{\Phi'\left(\xi\right)+\sqrt{\left[h\left(\xi\right)-\Phi^{2}\left(\xi\right)\right]^{2}+\left[\Phi'\left(\xi\right)\right]^{2}}}{h\left(\xi\right)-\Phi^{2}\left(\xi\right)}\right],
    \end{equation}
    \end{subequations}
then for two arbitrary constants $A$ and $B$
\begin{equation}\label{casesS0}
    u(x)=A u_{\mathtt{top}}(x) + B u_{\mathtt{bot}}(x),
\end{equation}
is a general solution of (\ref{Eq.main}) on $\left]x_-,x_+\right[$ if  and only if $\Phi(x)$ is a solution of (\ref{ELexp}) on this interval.    
\end{thmO}
\begin{corO}[Corollary 4 from \cite{Rudnicki}]\label{Cor2}\begin{subequations}\label{qmp}
    Functions 
    \begin{equation}\label{qmp1}
        \Theta_{\mathtt{top}}(x)=\frac{u'_{\mathtt{top}}(x)}{u_{\mathtt{top}}(x)}=\Phi\left(x\right)\frac{\Phi'\left(x\right)-\sqrt{\left[h\left(x\right)-\Phi^{2}\left(x\right)\right]^{2}+\left[\Phi'\left(x\right)\right]^{2}}}{h\left(x\right)-\Phi^{2}\left(x\right)},
    \end{equation}
      \begin{equation}\label{qmp2}
        \Theta_{\mathtt{bot}}(x)=\frac{u'_{\mathtt{bot}}(x)}{u_{\mathtt{bot}}(x)}=\Phi\left(x\right)\frac{\Phi'\left(x\right)+\sqrt{\left[h\left(x\right)-\Phi^{2}\left(x\right)\right]^{2}+\left[\Phi'\left(x\right)\right]^{2}}}{h\left(x\right)-\Phi^{2}\left(x\right)},
    \end{equation}\end{subequations}
    are two distinct solutions of a Ricatti equation
    \begin{equation}\label{Ricatti}
       \Theta'\left(x\right)+\Theta^2\left(x\right)+h(x)=0,
    \end{equation}
if  and only if $\Phi(x)$ is a solution of (\ref{ELexp}) on $\left]x_-,x_+\right[$. Moreover, the relations (\ref{specsol}) and (\ref{qmp}) can be uniquely inverted
 \begin{equation}\label{PhiSqrt}
   \Phi(x) = \sqrt{-\frac{u'_\mathtt{top}(x)u'_\mathtt{bot}(x)}{u_\mathtt{top}(x)u_\mathtt{bot}(x)}} =\sqrt{-\Theta_{\mathtt{top}}(x)\Theta_{\mathtt{bot}}(x)}. 
\end{equation}
\end{corO}

\section{Main results}\label{Sec2}
Our aim is to comprehensively generalize the former findings. As already listed, we are going to pursue this goal in three different, though mutually related, geometric ways which we call: (anti-)de Sitter, complex Riemannian and K{\"a}hler-Norden. To maintain readability we shall summarize the notation used in this manuscript.
\begin{notation}\label{notation}
    In \fullref{Table}{Table1}, I summarize notation used in \cite{Rudnicki} and in \fullref{Sec.}{subsechyp}, and list all variants of this notation introduced in \fullref{Sec.}{Sec2}. In the $2^\mathrm{nd}$ line I list sections in which the case under consideration is mainly being exposed. Most importantly, geodesic curves in explicit form, listed in $9^\mathrm{th}$ line, are denoted differently than the associated coordinates in the $3^\mathrm{rd}$ line. On the other hand, coordinates of the geodesic curves follow the latter notation, i.e. $\mathfrak{X}(\upsilon)$ is the second complex coordinate of the geodesic curve, even though $\textswab{X}\left(z\right)$ is the associated geodesic curve in explicit form. Finally, the bullet "$\bullet$" labeling solutions in the last line denotes all possible subscripts used, in particular, $\mathtt{top}$ and $\mathtt{bot}$.  
\end{notation}

\begin{table}
\begin{centering}
\caption{A summary of notation introduced in this manuscript. For details, see the notation remark at the beginning of \fullref{Sec.}{Sec2}.}
\begin{tabular}{|c|c|c|c|c|}
\hline \label{Table1}
Geometry & hyperbolic & (anti-)de Sitter  & \CWF{complex \\ Riemannian} & K{\"a}hler-Norden\tabularnewline
\hline 
\hline 
\CWF{Sections \\ \;}& \ref{subsechyp} \& \ref{Sec5} &\ref{subsec2.1}  \& \ref{proofs2.1} &\ref{subsec2.2}  \& \ref{proofs2.2} &\ref{subsec2.3}  \& \ref{proofs2.3}\tabularnewline
\hline
\CWF{Coordinates \\ \;}& $\left(x,\varPhi\right)$ & $\left(x,\varPsi\right)$ & $\left(z,\mathfrak{X}\right)$ & $\left(x,\varPhi,y,\varPsi\right)$\tabularnewline
\hline 
\CWF{Manifold\\ \;} & $\mathcal{H}_{0}$ & $\widetilde{\mathcal{H}}_{0}$ & $\mathcal{H}_{\mathbb{C}}$ & $\mathcal{H}_{\mathrm{KN}}$\tabularnewline
\hline
\CWF{Metric\\ \;} & $g_{h}$ & $\widetilde{g}_{h}$ & $G_{h}$ & $\hat{\mathbf{g}}_{h}$\tabularnewline
\hline 
\CWF{(Complex) \\ Riemannian\\ manifold} & $\mathbb{M}_{h}$ & $\widetilde{\mathbb{M}}_{h}^{\pm}$ & $\mathbb{M}_{h}^{\mathbb{C}}$ & $\mathbb{M}_{h}^{\mathrm{KN}}$\tabularnewline
\hline 
\CWF{Affine\\ parameter}& $s\in\mathbb{R}$ & $s\in\mathbb{R}$ & $\upsilon\in\mathbb{C}$ & $s\in\mathbb{R}$\tabularnewline
\hline 
\CWF{Christoffel\\ symbols} & $\Gamma$ & $\widetilde{\Gamma}$ & $\Upsilon$ & $\hat{\Gamma}$\tabularnewline
\hline 
\CWF{Geodesics in\\ explicit form} & $\Phi\left(x\right)$ & $\Psi\left(x\right)$ & $\textswab{X}\left(z\right)$ & $\Phi\left(x,y\right)$, $\Psi\left(x,y\right)$\tabularnewline
\hline 
\CWF{Solutions of\\ (\ref{Eq.mainC}) or (\ref{Eq.main})}& $u_{\bullet}\left(x\right)$ & $\widetilde{u}_{\bullet}\left(x\right)$ & $u_{\bullet}\left(z\right)$ & ---\tabularnewline
\hline 
\end{tabular}
\par\end{centering}
\end{table}

\subsection{Lorentzian (anti-)de Sitter geometry in two dimensions} \label{subsec2.1}
Upon a transformation $\varPhi\mapsto i \varPsi$, with $\varPsi>0$, the Riemannian metric $g_h$ becomes Lorentzian and degenerate when $\varPsi^2=-h(x)$. This brings a counterpart of \fullref{Definition}{DefM}.
\begin{defn}[Linear-$2^{nd}$-order-ODE (anti-)de Sitter geometry] \label{DefML}
Let $\widetilde{\mathcal{H}}_0 =\left\{\left(x,\varPsi\right)\in\mathbb{R}\times\mathbb{R}_{>0}\,|\,\varPsi^2\neq -h(x)\right\}$. Let
\begin{equation}\label{metricDS}
    \widetilde{g}_h = \frac{-\left(h(x)+\varPsi^2\right)^2dx^2+d\varPsi^2}{\varPsi^2},
\end{equation}
be a metric on $\widetilde{\mathcal{H}}_0$. Pairs $\widetilde{\mathbb{M}}_h^+ = \left(\widetilde{\mathcal{H}}_0, \widetilde{g}_h\right)$ and $\widetilde{\mathbb{M}}_h^- = \left(\widetilde{\mathcal{H}}_0, -\widetilde{g}_h\right)$, which only differ by a sign choice in front of $\widetilde{g}_h$, are two Lorentzian manifolds.
\end{defn}
Properties of $\widetilde{\mathbb{M}}_h^+$ and $\widetilde{\mathbb{M}}_h^-$, while vastly inherited from $\mathbb{M}_h$, offer a few interesting twists. 
\begin{lem} \label{LemL1} $\widetilde{\mathbb{M}}_h^+$ and $\widetilde{\mathbb{M}}_h^-$ have the following properties:
\begin{itemize} 
    \item[(i)] Sectional curvature of  $\widetilde{g}_h$ is equal to $-1$.  Thus, $\widetilde{\mathbb{M}}_h^+$ is locally a two dimensional anti-de Sitter space, while the manifold $\widetilde{\mathbb{M}}_h^-$ is locally a two dimensional de Sitter space.
    \item[(ii)] Both $\widetilde{\mathbb{M}}_h^+$ and $\widetilde{\mathbb{M}}_h^-$ share all geodesic curves.
    \item[(iii)] Geodesic curves on $\widetilde{\mathbb{M}}_h^+$ and $\widetilde{\mathbb{M}}_h^-$ in explicit form $\Psi(x)$ obey
    \begin{equation}\label{Thetageo}
\Psi''\left(x\right)=\frac{3\Psi^{2}\left(x\right)-h\left(x\right)}{\Psi^{2}\left(x\right)+h\left(x\right)}\frac{\left[\Psi'\left(x\right)\right]^{2}}{\Psi\left(x\right)}+\frac{h'\left(x\right)\Psi'\left(x\right)}{\Psi^{2}\left(x\right)+h\left(x\right)}-\frac{\Psi^{4}\left(x\right)-h^{2}\left(x\right)}{\Psi\left(x\right)}.
    \end{equation}
\end{itemize}
\end{lem}
The above results offer a justification for the term \textit{(anti-)de Sitter}, which aims to emphasize that in a two dimensional scenario the same or parallel results simultaneously hold in de Sitter, as well as in anti-de Sitter space. 

In an analogous way, we transfer \fullref{Theorem}{thm-main2} and \fullref{Corollary}{Cor2} to the Lorentzian geometry of $\widetilde{\mathbb{M}}_h^+$. By virtue of (ii) in \fullref{Lemma}{LemL1} we do not separately consider $\widetilde{\mathbb{M}}_h^-$ in the context of Eq. (\ref{Eq.main}).

\begin{thm}
Solutions of all linear second order ordinary differential equations (\ref{Eq.main}), with any $h(x)$ of class $C^1$, can be locally expressed through geodesic curves in  two dimensional anti-de Sitter geometry.  \label{thm-main2L}
    
More precisely, let  $h(x)$ be of class $C^1$, let $x_0\in\mathbb{R}$ and let 
$\Psi(x)$ be of class $C^2$ on some interval $\left]x_-,x_+\right[$, such that $x_-<x_0<x_+$ and $h\!\left(x\right)\neq-\Psi^{2}\!\left(x\right)$ on this interval. If $h\left(x\right)+\Psi^{2}\left(x\right)\not\equiv\pm\Psi'\left(x\right)$ on $\left]x_-,x_+\right[$ and
\begin{subequations}\label{specsolL}
\begin{equation}\label{specsol1L}
        \widetilde{u}_\mathtt{top}(x)=\exp\left[-\int_{x_0}^{x}\!\!d\xi\,\Psi\left(\xi\right)\frac{\Psi'\left(\xi\right)+i\sqrt{\left[h\left(\xi\right)+\Psi^{2}\left(\xi\right)\right]^{2}-\left[\Psi'\left(\xi\right)\right]^{2}}}{h\left(\xi\right)+\Psi^{2}\left(\xi\right)}\right],
    \end{equation}
    \begin{equation}\label{specsol2L}
        \widetilde{u}_\mathtt{bot}(x)=\exp\left[-\int_{x_0}^{x}\!\!d\xi\,\Psi\left(\xi\right)\frac{\Psi'\left(\xi\right)-i\sqrt{\left[h\left(\xi\right)+\Psi^{2}\left(\xi\right)\right]^{2}-\left[\Psi'\left(\xi\right)\right]^{2}}}{h\left(\xi\right)+\Psi^{2}\left(\xi\right)}\right],
    \end{equation}
    \end{subequations}
then for two arbitrary constants $A$ and $B$
\begin{equation}\label{casesS0L}
    u(x)=A \widetilde{u}_{\mathtt{top}}(x) + B \widetilde{u}_{\mathtt{bot}}(x),
\end{equation}
is a general solution of (\ref{Eq.main}) on $\left]x_-,x_+\right[$ if  and only if $\Psi(x)$ is a solution of (\ref{Thetageo}) on this interval.    
\end{thm}
\begin{cor}\label{Cor2L}\begin{subequations}\label{qmpL}
    Complex functions 
    \begin{equation}\label{qmp1L}
        \widetilde{\Theta}_{\mathtt{top}}(x)=\frac{\widetilde{u}'_{\mathtt{top}}(x)}{\widetilde{u}_{\mathtt{top}}(x)}=-\Psi\left(x\right)\frac{\Psi'\left(x\right)+i\sqrt{\left[h\left(x\right)+\Psi^{2}\left(x\right)\right]^{2}-\left[\Psi'\left(x\right)\right]^{2}}}{h\left(x\right)+\Psi^{2}\left(x\right)},
    \end{equation}
      \begin{equation}\label{qmp2L}
         \widetilde{\Theta}_{\mathtt{bot}}(x)=\frac{\widetilde{u}'_{\mathtt{bot}}(x)}{\widetilde{u}_{\mathtt{bot}}(x)}=-\Psi\left(x\right)\frac{\Psi'\left(x\right)-i\sqrt{\left[h\left(x\right)+\Psi^{2}\left(x\right)\right]^{2}-\left[\Psi'\left(x\right)\right]^{2}}}{h\left(x\right)+\Psi^{2}\left(x\right)},
    \end{equation}\end{subequations}
    are solutions of the Ricatti equation (\ref{Ricatti}) if  and only if $\Psi(x)$ is a solution of (\ref{Thetageo}) on $\left]x_-,x_+\right[$. Moreover, the relations (\ref{specsolL}) and (\ref{qmpL}) can be uniquely inverted
 \begin{equation}\label{PhiSqrtL}
   \Psi(x) = \sqrt{\frac{\widetilde{u}'_\mathtt{top}(x)\widetilde{u}'_\mathtt{bot}(x)}{\widetilde{u}_\mathtt{top}(x)\widetilde{u}_\mathtt{bot}(x)}} =\sqrt{ \widetilde{\Theta}_{\mathtt{top}}(x) \widetilde{\Theta}_{\mathtt{bot}}(x)}. 
\end{equation}
\end{cor}

A new aspect in \fullref{Theorem}{thm-main2L}, in comparison with \fullref{Theorem}{thm-main2}, is that now we have to include a requirement $h\left(x\right)+\Psi^{2}\left(x\right)\not\equiv\pm\Psi'\left(x\right)$, absent in  \fullref{Theorem}{thm-main2}. This also means that the solutions (\ref{qmp1L}) and (\ref{qmp2L}) coincide when $\left[h\left(x\right)+\Psi^{2}\left(x\right)\right]^{2}$ is  on $\left]x_-,x_+\right[$ identically equal to $\left[\Psi'\left(x\right)\right]^{2}$.

On the other hand, the subscripts "$\mathtt{top}/\mathtt{bot}$" in \fullref{Theorem}{thm-main2} refer to monotonicity of the solutions (\ref{specsol}). In \fullref{Theorem}{thm-main2L} this classification holds only when $\left|h\left(x\right)+\Psi^{2}\left(x\right)\right|<\left|\Psi'\left(x\right)\right|$. In this case, if $\Psi^2(x)$ lays on top/bottom of $-h(x)$, the function $\widetilde{u}_{\mathtt{top}}(x)$ is real and increasing/decreasing. In the opposite regime, the solutions (\ref{specsolL}) are oscillatory and mutually complex conjugate. An example is offered by the harmonic oscillator $h(x)=+\omega^2$ with $\omega>0$. The function  $\Psi_+(x)=\omega$ is a solution of (\ref{Thetageo}) which  for $x_0=0$ gives 
\begin{equation}
    \widetilde{u}_{\mathtt{top},+}(x)=e^{-i \omega x}, \qquad \widetilde{u}_{\mathtt{bot},+}(x)=e^{i \omega x}.
\end{equation}

We are only left with the "singular" case $h\left(x\right)+\Psi^{2}\left(x\right)=\pm\Psi'\left(x\right)$ for all $x\in\left]x_-,x_+\right[$, what means that either $\Psi\left(x\right)$ or $-\Psi\left(x\right)$ itself obeys the Ricatti equation (\ref{Ricatti}). Since, by assumption, $h\left(x\right)+\Psi^{2}\left(x\right)\neq 0$, these two possibilities are mutually exclusive. We then get a new feature of the discussed geometry.
\begin{thm} \label{LemL2}
    Solutions of the Ricatti equation (\ref{Ricatti}) with any $h(x)$ of class $C^1$ are geodesic curves in two dimensional anti-de Sitter geometry. More precisely, if $h\!\left(x\right)\neq-\Psi^{2}\!\left(x\right)$ on $\left]x_-,x_+\right[$ and either $\Psi(x)$ or $-\Psi(x)$ is a solution of (\ref{Ricatti}), then $\Psi(x)$ is a solution of (\ref{Thetageo}).
\end{thm}
Obviously, there is no equivalence between solutions of (\ref{Ricatti}) and solutions of (\ref{Thetageo}), as otherwise \fullref{Theorem}{thm-main2L} would be void. Let me also stress that the results of the whole \fullref{Sec.}{subsec2.1} describe an \textit{alternative} (and \textit{not} a \textit{complementary}) geometric representation of the solutions of Eq. (\ref{Eq.main}).

\subsection{Complex Riemannian geometry} \label{subsec2.2}
A pair  $\left(\mathcal{M},G\right)$ is a complex Riemannian manifold \cite{LeBrun1, LeBrun2}, also called a holomorphic Riemannian manifold \cite{ganchev}, if $\mathcal{M}$ is a complex manifold of some (complex) dimension $N$ and $G$ is a holomorphic metric on $\mathcal{M}$ \cite{LeBrun2}. In local holomorphic coordinates $\left(Z^1,\ldots, Z^N\right)$ such a metric is symmetric, non-degenerate and does not depend on the conjugate coordinates $\bar{Z}^a$, for all $a=1,\ldots,N$. The holomorphic metric $G$ admits a unique symmetric holomorphic affine connection locally expressed by Christoffel symbols $\Upsilon_{\:\:\:jk}^{i}\left(Z\right)$ derived from  $G$ in a standard way \cite{LeBrun1}. Moreover, all other notions and tools are directly transferred from (pseudo)-Riemannian geometry, bearing in mind the assumption of holomorphicity \cite{LeBrun1}. In particular, complex geodesic curves parametrized by a complex affine parameter $\upsilon$ (now the dot represents derivatives with respect $\upsilon$) also obey  \cite{LeBrun2} the standard geodesic equation $\ddot{Z}^i+\Upsilon_{\:\:\:jk}^{i}\left(Z\right)\dot{Z}^j\dot{Z}^k=0$. Finally, again in a standard way, one defines a holomorphic sectional curvature \cite{ganchev}. As we can see, complex Riemannian geometry offers a room naturally suitable for smooth and transparent extension of our former results to the full complex case. To this end we first adapt \fullref{Definition}{DefM}.
\begin{defn}[Linear-$2^{nd}$-order-ODE complex Riemannian geometry]\label{DefMC} Let $\mathcal{H}_\mathbb{C} =\left\{\left(z,\mathfrak{X}\right)\in\mathcal{U}\times\mathbb{C}\setminus\left\{0\right\}\,|\,\mathfrak{X}^2\neq h(z)\right\}$, where an open set $\mathcal{U}\subset\mathbb{C}$ and $h(z)$ is holomorphic on  $\mathcal{U}$. Let
\begin{equation}\label{metricC}
    G_h=\frac{\left(h(z)-\mathfrak{X}^2\right)^2dz^2+d\mathfrak{X}^2}{\mathfrak{X}^2},
\end{equation}
be a metric on $\mathcal{H}_\mathbb{C}$. A pair $\mathbb{M}_h^\mathbb{C}  = \left(\mathcal{H}_\mathbb{C} , G_h\right)$ is a complex Riemannian manifold.
\end{defn}
Analogously to \fullref{Lemma}{lem-hyp} we can prove the following.
\begin{lem}
The metric $G_h$ in the whole domain $\mathcal{H}_\mathbb{C}$ has constant holomorphic sectional curvature equal to $-1$. Therefore, it is locally holomorphical isometric to a complex 2-sphere of (complex) radius equal to imaginary unit $i$. \label{lem-hypC}
\end{lem}
The complex $N$-sphere of radius $ r\in\mathbb{C}\setminus\left\{0\right\}$ obeys
\begin{equation}
    \left(Z^1\right)^2+\ldots+\left(Z^N\right)^2+\left(Z^{N+1}\right)^2=r^2.
\end{equation}
It is a  holomorphic hypersurface \cite{ganchev}  of a complex Euclidean space of dimension $N+1$ parametrized by the coordinates $\left(Z^1,\ldots, Z^{N+1}\right)$.
\begin{rem}
    Similarly to the (anti)-de Sitter case, since the notion of positive definiteness does not apply to $G_h$, one can as well consider $\beta G_h$ with $\beta\in\mathbb{C}$. This scaling parameter does not influence geodesic curves. On the other hand, the holomorphic sectional curvature of $\beta G_h$ is $-\beta^{-1}$. Such geometry is locally holomorphical isometric to a complex 2-sphere of radius $i \sqrt{\beta}$. However, there is no qualitative influence of $\beta$ on the presented results.
\end{rem}

In the second step we also need a slight adjustment of \fullref{Definition}{defGEO}.
\begin{defn}\label{defGEOC} Let a geodesic curve on $\mathbb{M}_h^\mathbb{C}$ be parametrized by a complex affine parameter $\upsilon$ as  $\mathbb{C}\ni \upsilon \mapsto \left(z(\upsilon),\mathfrak{X}(\upsilon)\right)\in \mathcal{H}_\mathbb{C}$. If there exist open sets $\mathcal{V}_0\subset\mathbb{C}$ and $\mathcal{U}_0\subset\mathcal{U}$ such that $0\in\mathcal{V}_0$ and $z_0=z(0)\in\mathcal{U}_0$, and the map $\mathcal{V}_0\ni\upsilon\mapsto z(\upsilon)\in\mathcal{U}_0$ is such that $\dot{z}(\upsilon)\neq0$, then we  substitute $\mathfrak{X}(\upsilon)=\textswab{X}\left(z(\upsilon)\right)$, where $\textswab{X}\left(z\right)$  is a geodesic curve on $\mathbb{M}_h^\mathbb{C}$ represented in explicit form. 
\end{defn}
 
\begin{rem}\label{RemI4}
    Since the geodesic curve is subject to generic initial conditions, similarly as pointed out in \fullref{Remark}{RemI1}, we can always assure that $\dot{z}\left(0\right)=\dot{z}_0\neq 0$, thus,  geodesic curves  in explicit form can be defined at least locally. Moreover, as phrased in \fullref{Remark}{RemI2}, univalence of $z(\upsilon)$ is not a prerequisite to define $\textswab{X}\left(z\right)$. For the latter purpose it suffices that $\dot{z}(\upsilon)\neq 0$ for all $\upsilon\in\mathcal{V}_0$. For a single real variable both criteria are equivalent, while currently they are not (univalence is stricter). However, if $\dot{z}_0\neq 0$, then  due to inverse mapping theorem for holomorphic functions \cite{Kaup}, there also exist sets $\mathcal{V}_0$ and $\mathcal{U}_0$ for which $z(\upsilon)$ is biholomorphic.
\end{rem}

Moreover, geodesic curves are holomorphic because the connection is holomorphic by assumption. We automatically get a variant of \fullref{Proposition}{Prop1}.
\begin{prop}
Let $h(z)$ be holomorphic on $\mathcal{U}_0$. A geodesic curve on $\mathbb{M}_h^\mathbb{C}$ in explicit form $\textswab{X}(z)$ obeys \label{Prop1C}
\begin{equation}\label{ELexpC}
\textswab{X}''\left(z\right) =\frac{3\textswab{X}^{2}\left(z\right)+h\left(z\right)}{\textswab{X}^{2}\left(z\right)-h\left(z\right)}\frac{\left[\textswab{X}'\left(z\right)\right]^{2}}{\textswab{X}\left(z\right)}-\frac{h'\left(z\right)\textswab{X}'\left(z\right)}{\textswab{X}^{2}\left(z\right)-h\left(z\right)}+\frac{\textswab{X}^{4}\left(z\right)-h^{2}\left(z\right)}{\textswab{X}\left(z\right)}.
\end{equation}
\end{prop}

Bearing in mind a modified structure of solutions in (anti)-de Sitter geometry, which also applies now because $\left[h\left(z\right)-\textswab{X}^{2}\left(z\right)\right]^{2}+\left[\textswab{X}'\left(z\right)\right]^{2}$ might in principle vanish, we can extend former theorems to Eq. (\ref{Eq.mainC}). For consistency, we shall keep labeling based on the subscripts "$\mathtt{top}/\mathtt{bot}$" even though they reflect the original properties of the solutions only in special cases.

\begin{thm} Solutions of all linear second order ordinary differential equations (\ref{Eq.mainC}), with any $h(z)$ holomorphic on an open set $\mathcal{U}\subset\mathbb{C}$, can be locally expressed through geodesic curves in two dimensional geometry of a complex sphere with imaginary unit radius.  \label{thm-main2C}
    
More precisely, let  $h(z)$ be holomorphic on $\mathcal{U}$, let $z_0\in\mathcal{U}_0\subset\mathcal{U}$, let 
$\textswab{X}(z)$ be holomorphic on $\mathcal{U}_0$ and let $h\!\left(z\right)\neq\textswab{X}^{2}\!\left(z\right)$ on $\mathcal{U}_0$. Let $\gamma_z$ be an arbitrary holomorphic curve in $\mathcal{U}_0$ which starts at $z_0$ and ends at $z$. If $\left[h\left(z\right)-\textswab{X}^{2}\left(z\right)\right]^{2}+\left[\textswab{X}'\left(z\right)\right]^{2}\not\equiv0$ on $\mathcal{U}_0$ and
\begin{subequations}\label{specsolC}
\begin{equation}\label{specsol1C}
        u_\mathtt{top}(z)=\exp\left[\int_{\gamma_z}\!\!d\zeta\,\textswab{X}\left(\zeta\right)\frac{\textswab{X}'\left(\zeta\right)-\sqrt{\left[h\left(\zeta\right)-\textswab{X}^{2}\left(\zeta\right)\right]^{2}+\left[\textswab{X}'\left(\zeta\right)\right]^{2}}}{h\left(\zeta\right)-\textswab{X}^{2}\left(\zeta\right)}\right],
    \end{equation}
    \begin{equation}\label{specsol2C}
        u_\mathtt{bot}(z)=\exp\left[\int_{\gamma_z}\!\!d\zeta\,\textswab{X}\left(\zeta\right)\frac{\textswab{X}'\left(\zeta\right)+\sqrt{\left[h\left(\zeta\right)-\textswab{X}^{2}\left(\zeta\right)\right]^{2}+\left[\textswab{X}'\left(\zeta\right)\right]^{2}}}{h\left(\zeta\right)-\textswab{X}^{2}\left(\zeta\right)}\right],
    \end{equation}
    \end{subequations}
then for two arbitrary constants $A$ and $B$
\begin{equation}\label{casesS0C}
    u(z)=A u_{\mathtt{top}}(z) + B u_{\mathtt{bot}}(z),
\end{equation}
is a general solution of (\ref{Eq.mainC}) on $\mathcal{U}_0$ if  and only if $\textswab{X}(z)$ is a solution of (\ref{ELexpC}) on this open set.  The functions (\ref{specsolC}) do not depend on the shape of the curve $\gamma_z$.  
\end{thm}
\begin{cor}\label{Cor2C}\begin{subequations}\label{qmpC}
    Complex functions 
    \begin{equation}\label{qmp1C}
        \Theta_{\mathtt{top}}(z)=\frac{u'_{\mathtt{top}}(z)}{u_{\mathtt{top}}(z)}=\textswab{X}\left(z\right)\frac{\textswab{X}'\left(z\right)-\sqrt{\left[h\left(z\right)-\textswab{X}^{2}\left(z\right)\right]^{2}+\left[\textswab{X}'\left(z\right)\right]^{2}}}{h\left(z\right)-\textswab{X}^{2}\left(z\right)},
    \end{equation}
      \begin{equation}\label{qmp2C}
        \Theta_{\mathtt{bot}}(z)=\frac{u'_{\mathtt{bot}}(z)}{u_{\mathtt{bot}}(z)}=\textswab{X}\left(z\right)\frac{\textswab{X}'\left(z\right)+\sqrt{\left[h\left(z\right)-\textswab{X}^{2}\left(z\right)\right]^{2}+\left[\textswab{X}'\left(z\right)\right]^{2}}}{h\left(z\right)-\textswab{X}^{2}\left(z\right)},
    \end{equation}\end{subequations}
    are two solutions of the complex Ricatti equation
    \begin{equation}\label{RicattiC}
       \Theta'\left(z\right)+\Theta^2\left(z\right)+h(z)=0,
    \end{equation}
if  and only if $\textswab{X}(z)$ is a solution of (\ref{ELexpC}) on $\mathcal{U}_0$. Moreover, the relations (\ref{specsolC}) and (\ref{qmpC}) can be uniquely inverted
 \begin{equation}\label{PhiSqrtC}
   \textswab{X}(z) = \sqrt{-\frac{u'_\mathtt{top}(z)u'_\mathtt{bot}(z)}{u_\mathtt{top}(z)u_\mathtt{bot}(z)}} =\sqrt{-\Theta_{\mathtt{top}}(z)\Theta_{\mathtt{bot}}(z)}. 
\end{equation}
\end{cor}
\begin{thm} \label{LemL2C}
    Solutions of the complex Ricatti equation (\ref{RicattiC}) with any holomorphic $h(z)$ are geodesic curves in two dimensional geometry of a complex sphere with imaginary unit radius. More precisely, if $h\!\left(z\right)\neq\textswab{X}^{2}\!\left(z\right)$ on $\mathcal{U}_0$ and either $i\textswab{X}(z)$ or $-i\textswab{X}(z)$ is a solution of (\ref{RicattiC}), then $\textswab{X}(z)$ is a solution of (\ref{ELexpC}).
\end{thm}

Before moving to the last part, let us discuss in more detail a construction of the curve $\gamma_z$. This will turn out to be essential for the final statements of \fullref{Sec.}{subsec2.3}. Every such path in $\mathcal{U}_0$ is an image of a path $\gamma^\circ_z$ in $\mathcal{V}_0$. Let us introduce a real affine parameter $s\in\left[0,1\right]$ and parametrize $\gamma^\circ_z$ as $\upsilon(s)\in\mathcal{V}_0$, with $\upsilon(0)=0$ and $\upsilon(1)=\upsilon_z$. Without loss of generality we can always rescale the holomorphic affine parameter $\upsilon$ such that $\upsilon_z=1$.

If, as in (\ref{specsolC}), we work with a complex variable $\zeta\in\mathcal{U}_0$, taken to be a holomorphic function $\zeta(\upsilon)$, we further get $\zeta(0)=z_0$ and $\zeta(1)=z$. In this way, we construct a curve $\left[0,1\right]\rightarrow\mathcal{U}_0$ of the form $\zeta\left(\upsilon(s)\right)$. As a consequence, we can rewrite
\begin{equation}\label{integralGamma}
    \int_{\gamma_z}\!\!d\zeta\, \Theta_{\mathtt{top}/\mathtt{bot}}(\zeta)=\int_0^1 \!\!ds\,\frac{d\zeta\left(\upsilon(s)\right)}{d\upsilon}\frac{d\upsilon(s)}{ds}\Theta_{\mathtt{top}/\mathtt{bot}}\left(\zeta\left(\upsilon(s)\right)\right).
\end{equation}
As already mentioned, the above integral which enters (\ref{specsolC}) does not depend on the shape of $\gamma_z$ as long as $\zeta(\upsilon)$ is holomorphic. More importantly, the integrand depends on $\textswab{X}(\zeta(\upsilon))$ and its derivative, taken at $\upsilon=\upsilon(s)$. Therefore, in relation to \fullref{Remark}{RemI4}, we recognize that the sole requirement of non-vanishing derivative of $z(\upsilon)$ makes the formulas (\ref{specsolC}) meaningful, even if the function $z(\upsilon)$ is not biholomorphic. This follows from the fact that one can take $\gamma_z$ to be a first coordinate of the geodesic curve in question.

\subsection{Pseudo-Riemannian K{\"a}hler-Norden geometry in four dimensions} \label{subsec2.3} In the last part we link $N$-dimensional complex Riemannian geometry with its underlying pseudo-Riemannian real geometry, classified by Berger \cite{Berger1955} according to its holonomy group $SO\left(N,\mathbb{C}\right)$ \cite{Bryant}. Given local holomorphic coordinates $\left(Z^1,\ldots, Z^N\right)$ on $\mathcal{M}$, we split them into their real and imaginary parts $Z^k=x^k+i y^k$ and define \cite{Borowiec} \begin{equation}
X^{k}=\begin{cases}
x^{k} & k=1,\ldots,N\\
y^{k-N} & k=N+1,\ldots,2N
\end{cases},
\end{equation}
to be $2N$ real local coordinates on $\mathcal{M}$. Every complex Riemannian manifold  also admits a real pseudo-Riemannian K{\"a}hler-Norden metric $\hat{\mathbf{g}}$ defined by the formula \cite{Borowiec} (see also \cite{KN1,KN2,KN3})
\begin{equation}\label{correpon}
\hat{\mathbf{g}}=\frac{\alpha}{2}\left(G+\bar{G}\right),\qquad\hat{\mathbf{g}}_{kl}dX^{k}dX^{l}=\alpha\textrm{Re}\left[G_{ab}dZ^{a}dZ^{b}\right],
\end{equation}
where $\alpha$ is a real scaling constant (in \cite{Borowiec} set as $\alpha=2$) and $\bar{G}=\bar{G}_{ab}d\bar{Z}^{a}d\bar{Z}^{b}$. The metric $\hat{\mathbf{g}}$ is balanced, i.e. of signature $\left(-,-,+,+\right)$. While we call the pair $\left(\mathcal{M},G\right)$ as complex Riemannian manifold, we shall call the $\left(\mathcal{M},\hat{\mathbf{g}}\right)$ pair a K{\"a}hler-Norden manifold. The latter also formally requires an almost complex structure. Importantly, from an interpretational point of view, this manifold is anti-K{\"a}hlerian, i.e. the metric $\hat{\mathbf{g}}$, also called a "B-metric" \cite{Bmetric}, is anti-Hermitian with respect to the almost complex structure. This is consistent with findings from \fullref{Sec.}{subsec2.2} where, even though not spelled out explicitly, it was not possible to employ a model of a complex hyperbolic geometry \cite{Comhyp} which requires a hermitian metric. We need to introduce a real four dimensional extension of \fullref{Definition}{DefM}, setting $\alpha=1$ in Eq. (\ref{correpon}).
\begin{defn}[Linear-$2^{nd}$-order-ODE  K{\"a}hler-Norden geometry] Let 
\begin{equation}
    \mathcal{H}_\mathrm{KN} =\left\{\left(x,\varPhi,y,\varPsi\right)\in\mathbb{R}^4\,|\,\left(\varPhi+i\varPsi\neq 0\right)\wedge
   \left(h\left(x,y\right)\neq\left(\varPhi+i \varPsi\right)^{2}\right)\right\},
\end{equation}
where $h\left(x,y\right)=h_{\Re}\left(x,y\right)+i h_{\Im}\left(x,y\right)$, and $h_{\Re}\left(x,y\right)$ together with $h_{\Im}\left(x,y\right)$ are defined on an open set $\mathcal{U}\subset\mathbb{R}^2$ and fulfill
\begin{equation}\label{CauchyRiemann}
\frac{\partial h_{\Re}\left(x,y\right)}{\partial x}=\frac{\partial h_{\Im}\left(x,y\right)}{\partial y},\qquad\frac{\partial h_{\Re}\left(x,y\right)}{\partial y}=-\frac{\partial h_{\Im}\left(x,y\right)}{\partial x}.
\end{equation}
Let a metric on $\mathcal{H}_\mathrm{KN}$ be 
\begin{eqnarray}\label{metricKN}
\hat{\mathbf{g}}_{h} & = & \frac{\Delta_{-}\left[\Delta_{+}^{2}+h_{\Re}^{2}-h_{\Im}^{2}\right]+4\varPhi\varPsi h_{\Re}h_{\Im}-2\Delta_{+}^{2}h_{\Re}}{\Delta_{+}^{2}}\left(dx^{2}-dy^{2}\right)\nonumber\\
 & - & 4\frac{\left[\varPhi h_{\Im}-\varPsi\left(\Delta_{+}+h_{\Re}\right)\right]\left[\varPsi h_{\Im}-\varPhi\left(\Delta_{+}-h_{\Re}\right)\right]}{\Delta_{+}^{2}}dxdy\nonumber\\
 & + & \frac{\Delta_{-}\left(d\varPhi^{2}-d\varPsi^{2}\right)+4\varPhi\varPsi d\varPhi d\varPsi}{\Delta_{+}^{2}},
\end{eqnarray}
where $\Delta_{\pm}=\Delta_{\pm}\left(\varPhi,\varPsi\right)=\varPhi^{2}\pm\varPsi^{2}$
and $h_{\Re/\Im}=h_{\Re/\Im}\left(x,y\right)$. \label{DefMKN}
A pair $\mathbb{M}_h^\mathrm{KN}  = \left(\mathcal{H}_\mathrm{KN}, \hat{\mathbf{g}}_h\right)$ is a K{\"a}hler-Norden pseudo-Riemannian manifold.
\end{defn}
We notice that (\ref{CauchyRiemann}) are Cauchy-Riemann equations which assure that $h\left(x,y\right)$
is a holomorphic function of $z=x+iy$ on $\mathcal{U}$. We also observe that if $y=0$ (so that $z=x$) and $h_{\Im}\left(x,0\right)=0$, then  
\begin{itemize}
    \item If $\varPsi=0$,  then $\hat{\mathbf{g}}_{h} = g_h$ with $g_h$ from \fullref{Definition}{DefM}.
    \item If  $\varPhi=0$, then $\hat{\mathbf{g}}_{h} = \widetilde{g}_h$ with $\widetilde{g}_h$ from \fullref{Definition}{DefML}.
\end{itemize}
In this way $\hat{\mathbf{g}}_{h}$ subsumes two dimensional submanifolds, which were previously used to describe the real scenario. Therefore, \fullref{Definition}{DefMKN} is automatically a generalization of the former results. 

Even though, the metric (\ref{metricKN}) on its own in principle describes a very cumbersome geometry, due to the Cauchy-Riemann equations (\ref{CauchyRiemann}) it just reduces to a different description of geometry given by (\ref{metricC}).

Still, while complex sectional curvature of $G_h$ is constant, it does not imply that standard (real) sectional curvature of $\hat{\mathbf{g}}_{h}$ would have this property. Since it is not the case, its sectional curvature is actually necessarily unbounded \cite{Kulkarni}. However, the property of being an Einstein metric \cite{Besse}, defined as $R=\eta g$ (where $g$ is a metric, $R$ is its Ricci tensor, and $\eta$ is a real constant), is inherited.
\begin{lem}\label{Einstein}
    Metric $\hat{\mathbf{g}}_{h}$ is Einstein with a constant $\eta=-2$. Therefore, its Ricci scalar is equal to $-8$.
\end{lem}
We conclude our findings with the following theorem.
\begin{thm}\label{GeoInherited}
    Solutions of all linear second order ordinary differential equations (\ref{Eq.mainC}), with any $h(z)$ holomorphic on an open set $\mathcal{U}\subset\mathbb{C}$, can be locally expressed through geodesic curves in  K{\"a}hler-Norden geometry of a four dimensional pseudo-Riemannian Einstein manifold $\mathbb{M}_h^\mathrm{KN}$.  
    
More precisely, let $\mathbb{R}\ni s \mapsto \left(x(s),\varPhi(s),y(s),\varPsi(s)\right)\in \mathcal{H}_\mathrm{KN}$ be a geodesic curve with a real affine parameter $s$, let $\mathcal{V}_0$ be convex and let the path $\gamma_z^\circ$ in  $\mathcal{V}_0$ be defined as $\upsilon(s)= s$. Then, inside the solutions (\ref{specsolC}), which build a general solution (\ref{casesS0C}) of (\ref{Eq.mainC}), one can decompose holomorphic geodesic curves in explicit form as $\textswab{X}(z) = \Phi\left(x,y\right)+i \Psi\left(x,y\right)$, where:\begin{subequations}
\begin{equation}
    x(s)=\mathrm{Re}\left[z(s)\right],\qquad y(s)=\mathrm{Im}\left[z(s)\right], 
\end{equation}
\begin{equation}
   \Phi\left(x(s),y(s)\right)= \varPhi(s)=\mathrm{Re}\left[\mathfrak{X}(s)\right]=\mathrm{Re}\left[\textswab{X}\left(z(s)\right)\right],
\end{equation}
\begin{equation}
    \Psi\left(x(s),y(s)\right)=\varPsi(s)=\mathrm{Im}\left[\mathfrak{X}(s)\right]=\mathrm{Im}\left[\textswab{X}\left(z(s)\right)\right],
\end{equation}
\end{subequations}
and $\mathbb{R}\ni s \mapsto \left(z(s),\mathfrak{X}(s)\right)\in \mathcal{H}_\mathbb{C}$ is a geodesic curve on $\mathbb{M}_h^\mathbb{C}$ with the  holomorphic parameter $\upsilon$ restricted to the real axis.
\end{thm}

\section{Conclusions}\label{Sec3}
In \cite{Rudnicki} I established a mutual connection between real linear second order ordinary differential equations and two dimensional hyperbolic geometry. A foremost question left for the sequel manuscript was a generalization of the results to the complex scenario. In the current paper this goal has been achieved, provided that the function $h(x)$ is replaced by a holomorphic function $h(z)$. It turns out that the two dimensional hyperbolic geometry generalizes to the two dimensional complex sphere (see \fullref{Theorem}{thm-main2C}), described within complex Riemannian differential geometry (see \fullref{Lemma}{lem-hypC}). Such a manifold is characterized by a constant (complex) holomorphic sectional curvature, not necessarily equal to $1$. As already emphasized in Sec. 1 of \cite{Rudnicki}, such results are of a local and coordinate-dependent character, what follows from the very nature of the pursued approach. For example, in this particular case the manifold underpinning Eq. (\ref{Eq.mainC}) is "only" locally isometric to the complex sphere. In this way, we can see that further open questions in this research program could refer to a topological character of the  discovered geometric identifications. In addition, a related question not tackled here is about an explicit form of local diffeomorphisms between the (complex) Riemannian manifolds. Previously, we constructed such a diffeomorphism for the hyperbolic geometry of Eq. (\ref{Eq.main}). Now, we do not look for the generalization of Lemma 6 from  \cite{Rudnicki}, believing it shall be straightforward (though, still nontrivial).

For readers who prefer real differential geometry, we show that the main result is in "one to one" correspondence  (see \fullref{Theorem}{GeoInherited}) with pseudo-Riemannian four dimensional geometry of the so-called K{\"a}hler-Norden type. The family of metrics (\ref{metricKN}) can be considered a base metric for the developed theory which is equivalent to the complex Riemannian metric (\ref{metricC}), and reduces to both (\ref{metric}) and (\ref{metricDS}) which are relevant for the real equation (\ref{Eq.main}). It is worth stressing that we started with a complementary description of the real scenario  (see \fullref{Theorem}{thm-main2L}), which boils down to the (anti-)de Sitter geometry  (see \fullref{Lemma}{LemL1}).

A peculiar feature of (anti-)de Sitter and complex Riemannian geometry, absent in the case of the hyperbolic geometry, is that they admit a "singular" scenario in which the geodesic curves in explicit form are just solutions of the Ricatti equation (see \fullref{Theorem}{LemL2} and \fullref{Theorem}{LemL2C}).

I believe the results collected in the two manuscripts do not exhaust the topic they introduce, but rather, offering a comprehensive approach covering "regular" scenarios add to (almost) textbook knowledge in the field of ordinary differential equations. Especially in the complex case, there are several open questions for future consideration, for example, pertaining to situations in which $h(z)$ is not holomorphic. Moreover, the geometric features of the metric (\ref{metricKN}) remain largely unexplored, except results of \fullref{Lemma}{Einstein}.

\section{Proofs}
\label{Sec4}
\subsection{Proofs for (anti-)de Sitter geometry} \label{proofs2.1}
Technically speaking, all results in \fullref{Sec.}{subsec2.1} can be obtained from those of \fullref{Sec.}{subsechyp} by a replacement $\varPhi\mapsto i \varPsi$ treated the same way as a coordinate transformation in complex domain, as long as one remembers if an object under consideration is a scalar, a tensor or follows a different transformation pattern (like the Christoffel symbols do). Obviously, one can always perform an independent checkup starting from $\widetilde{g}_h$, but for the sake of brevity we are not following that route. Only in the proof of \fullref{Lemma}{LemL1} we report such a detailed calculation, which anyway is patterned after the proof of \fullref{Lemma}{lem-hyp} from \cite{Rudnicki}, just to avoid confusion in the transformation of Christoffel symbols.

\begin{proof}[Proof of \textup{\fullref{Lemma}{LemL1}}]
Given a metric $g$, its Christoffel symbols and Riemann tensor, both defined in a standard way, while described in some coordinates (Einstein summation convention applies) are as follows \cite{Besse}: 
\begin{subequations}
\begin{equation}
\Gamma_{\:\:\:jk}^{i}=\frac{1}{2}g^{il}\left(g_{lj,k}+g_{lk,j}-g_{jk,l}\right),
\end{equation}
\begin{equation}
R_{\:\:\:jkl}^{i}=\Gamma_{\:\:\:jl,k}^{i}-\Gamma_{\:\:\:jk,l}^{i}+\Gamma_{\:\:\:mk}^{i}\Gamma_{\:\:\:\:jl}^{m}-\Gamma_{\:\:\:ml}^{i}\Gamma_{\:\:\:\:jk}^{m}.
\end{equation}
\end{subequations}

In the chosen
coordinates $\left(x,\varPsi\right)$ the metric $\widetilde{g}_h$ is (I skip the label $h$)
\begin{subequations}
\begin{equation}
    \widetilde{g}_{xx} = -\frac{\left(h\left(x\right)+\varPsi^{2}\right)^{2}}{\varPsi^{2}},\quad \widetilde{g}_{\varPsi\varPsi}=\frac{1}{\varPsi^{2}},\quad \widetilde{g}_{x\varPsi}=\widetilde{g}_{\varPsi x}=0.
\end{equation}
Therefore:
\begin{equation}
\widetilde{\Gamma}_{\:\:\:xx}^{x}=\frac{h'\left(x\right)}{h\left(x\right)+\varPsi^{2}},\quad\widetilde{\Gamma}_{\:\:\:x\varPsi}^{x}=\widetilde{\Gamma}_{\:\:\:\varPsi x}^{x}=\frac{\varPsi^{2}-h\left(x\right)}{\varPsi\left[\varPsi^{2}+h\left(x\right)\right]},
\end{equation}
\begin{equation}
\widetilde{\Gamma}_{\:\:\:xx}^{\varPsi}=\frac{\varPsi^{4}-h^{2}\left(x\right)}{\varPsi},\quad\widetilde{\Gamma}_{\:\:\:\varPsi\varPsi}^{\varPsi}=-\frac{1}{\varPsi},
\end{equation}
\begin{equation}
\widetilde{R}_{\:\:\:\varPsi x\varPsi}^{x}=-\widetilde{R}_{\:\:\:\varPsi\varPsi x}^{x}=-\frac{1}{\varPsi^{2}},\quad \widetilde{R}_{\:\:\:xx\varPsi}^{\varPsi}=-\widetilde{R}_{\:\:\:x\varPsi x}^{\varPsi}=-\frac{\left(h\left(x\right)+\varPsi^{2}\right)^{2}}{\varPsi^{2}},
\end{equation}
\begin{equation}
\widetilde{R}_{x\varPsi x\varPsi}=-\widetilde{R}_{x\varPsi\varPsi x}=-\widetilde{R}_{\varPsi xx\varPsi}=\widetilde{R}_{\varPsi x\varPsi x}=\frac{\left(h\left(x\right)+\varPsi^{2}\right)^{2}}{\varPsi^{4}}.
\end{equation}
I have only listed non-vanishing coefficients. The Ricci tensor $\widetilde{R}_{ij}=\widetilde{R}_{\:\:\:ikj}^{k}$ is
\begin{equation}
\widetilde{R}_{xx}=\frac{\left(h\left(x\right)+\varPsi^{2}\right)^{2}}{\varPsi^{2}},\quad \widetilde{R}_{\varPsi\varPsi}=-\frac{1}{\varPsi^{2}},\quad \widetilde{R}_{x\varPsi}=\widetilde{R}_{\varPsi x}=0,
\end{equation}    
\end{subequations}
and turns out to be equal to $-\widetilde{g}_{h}$. 

Just for comparison we provide non-vanishing coefficients of the Christoffel symbols for $g_h$:\begin{subequations}\label{ChristoOrig}
\begin{equation}
\Gamma_{\:\:\:xx}^{x}=\frac{h'\left(x\right)}{h\left(x\right)-\varPhi^{2}},\quad\Gamma_{\:\:\:x\varPhi}^{x}=\Gamma_{\:\:\:\varPhi x}^{x}=\frac{\varPhi^{2}+h\left(x\right)}{\varPhi\left[\varPhi^{2}-h\left(x\right)\right]},
\end{equation}
\begin{equation}
\Gamma_{\:\:\:xx}^{\varPhi}=\frac{h^{2}\left(x\right)-\varPhi^{4}}{\varPhi},\quad\Gamma_{\:\:\:\varPhi\varPhi}^{\varPhi}=-\frac{1}{\varPhi},
\end{equation}
As before, all other coefficients vanish. With that background we are ready to prove assertions of \fullref{Lemma}{LemL1}.
\end{subequations}

\textit{(i)}. The value of the
sectional curvature of $\widetilde{g}_h$ is
\begin{equation}
\widetilde{K}=\frac{\widetilde{R}_{x\varPsi x\varPsi}}{\widetilde{g}_{xx}\widetilde{g}_{\varPsi\varPsi}-\widetilde{g}_{x\varPsi}\widetilde{g}_{\varPsi x}}\equiv-1.
\end{equation} 
Since it is a constant, $\widetilde{\mathbb{M}}^+_h$ is a locally maximally symmetric Lorentzian manifold. Moreover, as $\widetilde{K}$ is negative, we locally obtain anti-de Sitter space. Since, similarly to the case of $\mathbb{M}_h$, we know that the manifold in question is generally not complete, this property is not global. 

The manifold $\widetilde{\mathbb{M}}^-_h$ only differs from $\widetilde{\mathbb{M}}^+_h$ by a minus sign in front of the metric tensor. Obviously, the sectional curvature of $-\widetilde{g}_h$ is equal to $1$. This only difference transforms anti-de Sitter space into the de Sitter space. 

\textit{(ii)}. Since the Christoffel symbols are invariant under multiplication of a metric tensor by a constant, geodesic curves in both manifolds under consideration fulfill exactly the same geodesic equation.

\textit{(iii)}. Eq. (\ref{Thetageo}) follows from \fullref{Proposition}{Prop1} (its proof can be found in \fullref{Sec.}{Sec5}) in two equivalent ways. In the first approach, since the sets of vanishing Christoffel symbols for $\widetilde{g}_h$ and $g_h$ are the same, we immediately obtain (\ref{extradop}) with $\Phi(x)$ relabeled to $\Psi(x)$ and the Christoffel symbols $\Gamma$ replaced by $\widetilde{\Gamma}$ with the coordinate label $\varPhi$ replaced by $\varPsi$ (imaginary unit was already "absorbed" into new Christoffel symbols). Alternatively, we take  (\ref{ELexp}) and make the original replacement $\Phi(x)\mapsto i \Psi(x)$. Both approaches swiftly give the final result. 
\end{proof}

\begin{proof}[Proof of \textup{\fullref{Theorem}{thm-main2L}}]
Equations (\ref{specsolL}) follow immediately from (\ref{specsol}) by virtue of the replacement $\Phi(x)\mapsto i \Psi(x)$. Therefore, due to \fullref{Theorem}{thm-main2} (its proof can be found in \fullref{Sec.}{Sec5}), both (\ref{specsol1L}) and (\ref{specsol2L}) are particular solutions of (\ref{Eq.main}). These solutions are linearly independent if $h\left(x\right)+\Psi^{2}\left(x\right)\not\equiv\pm\Psi'\left(x\right)$, as assumed.
\end{proof}

\begin{proof}[Proof of \textup{\fullref{Corollary}{Cor2L}}]
Equations (\ref{qmpL}) and (\ref{PhiSqrtL}) follow immediately from (\ref{qmp}) and (\ref{PhiSqrt}) by virtue of the replacement $\Phi(x)\mapsto i \Psi(x)$. 
\end{proof}

\begin{proof}[Proof of \textup{\fullref{Theorem}{LemL2}}]
Let $\pm\Psi(x)$  be a solution of the Ricatti equation (\ref{Ricatti}). Since $h(x)\neq -\Psi^2(x)$, $\Psi'\left(x\right)$ never changes the sign and these possibilities are mutually exclusive, i.e. either $\Psi(x)$ or $-\Psi(x)$ solves (\ref{Ricatti}).

 On the other hand, under the above assumption, from (\ref{qmpL}) we find that 
\begin{equation}
     \widetilde{\Theta}_{\mathtt{top}}(x)= \widetilde{\Theta}_{\mathtt{bot}}(x)=\pm \Psi(x).
\end{equation}
Therefore, $\widetilde{\Theta}_{\mathtt{top}}(x)$ and $\widetilde{\Theta}_{\mathtt{bot}}(x)$ are as well solutions of (\ref{Ricatti}). According to \fullref{Corollary}{Cor2L}, $\Psi(x)$ is a solution of (\ref{Thetageo}).
\end{proof}

\subsection{Proofs for complex Riemannian geometry}\label{proofs2.2}
Since we only work with (locally) holomorphic functions, several results in \fullref{Sec.}{subsec2.2} can be obtained from those of \fullref{Sec.}{subsechyp} by a replacement $x \mapsto z$, $\varPhi \mapsto \mathfrak{X}$ and $\Phi(x) \mapsto \textswab{X}(z)$ accompanied by a careful account of differences brought by complex analysis.

\begin{proof}[Proof of \textup{\fullref{Lemma}{lem-hypC}}]
Geometry of $G_h$ is in local coordinates $\left(z,\mathfrak{X}\right)$ described by the very same formulas as provided in the proof of \fullref{Lemma}{lem-hyp}, to be found in \cite{Rudnicki} (see also proof of \fullref{Lemma}{LemL1} given above), with a minor replacement $x \mapsto z$ and $\varPhi \mapsto \mathfrak{X}$ (and remembering that $\mathfrak{X}\in\mathbb{C}\setminus\left\{0\right\}$). Just as an excerpt from these results (and for further use) I provide all non-vanishing Christoffel symbols for $G_h$:
\begin{subequations}\label{ChristoCom}
\begin{equation}
\Upsilon_{\:\:\:zz}^{z}=\frac{h'\left(z\right)}{h\left(z\right)-\mathfrak{X}^{2}},\quad\Upsilon_{\:\:\:z\mathfrak{X}}^{z}=\Upsilon_{\:\:\:\mathfrak{X} z}^{z}=\frac{\mathfrak{X}^{2}+h\left(z\right)}{\mathfrak{X}\left[\mathfrak{X}^{2}-h\left(z\right)\right]},
\end{equation}
\begin{equation}
\Upsilon_{\:\:\:zz}^{\mathfrak{X}}=\frac{h^{2}\left(z\right)-\mathfrak{X}^{4}}{\mathfrak{X}},\quad\Upsilon_{\:\:\:\mathfrak{X}\mathfrak{X}}^{\mathfrak{X}}=-\frac{1}{\mathfrak{X}}.
\end{equation}
\end{subequations}
Therefore, holomorphic sectional curvature of $G_h$, defined in Sec. 4 of \cite{ganchev}, must be equal to $-1$, by analogy with the case of $g_h$. In the same vein, the Ricci tensor of $G_h$ is equal to $-G_h$.

According to Proposition 4.1 from \cite{ganchev}, the holomorphic sectional curvature of the complex sphere with radius $r$ is equal to $r^{-2}$. For the choice $r=i$, such a sphere has the same holomorphic sectional curvature as $\mathbb{M}_h^\mathbb{C}$. Therefore, according to Theorem 5.3 from \cite{ganchev} both complex Riemannian manifolds in question are locally holomorphical isometric to each other.
\end{proof}

\begin{proof}[Proof of \textup{\fullref{Proposition}{Prop1C}}]
Since $h(z)$ is holomorphic,  Eq. (\ref{ELexpC}) follows immediately from \fullref{Proposition}{Prop1} (its proof can be found in \fullref{Sec.}{Sec5}) upon a replacement $x\mapsto z$ and $\Phi(x)\mapsto \textswab{X}(z)$.
\end{proof}

\begin{proof}[Proof of \textup{\fullref{Theorem}{thm-main2C}}]
Equations (\ref{specsolC}) are direct extensions of (\ref{specsol}) to the complex domain. The only technical difference lies in the range of integration (a curve in the complex plane instead of an interval on the real line), however, since the integrands in (\ref{specsolC}) given by (\ref{qmpC}) are holomorphic, every fixed curve $\gamma_z^1$ can be replaced by any other curve $\gamma_z^2$, as long as both live in $\mathcal{U}_0$. The difference would then be given by the integral over a loop $\gamma_z^1\cup\gamma_z^2$, which vanishes due to Cauchy theorem \cite{Kaup}. Thus, we already conclude that the functions  (\ref{specsolC}) do not depend on the shape of $\gamma_z$. In the same way we justify that (\ref{qmpC}) are just derivatives with respect to $z$ of the integrals from (\ref{specsolC}).

Therefore, the proof that  both (\ref{specsol1C}) and (\ref{specsol2C}) are particular solutions of (\ref{Eq.mainC}) is literally the same as the proof of \fullref{Theorem}{thm-main2} (which can be found in \fullref{Sec.}{Sec5}). One just needs to replace  $\Phi(x)\mapsto \textswab{X}(z)$. These solutions are linearly independent if $\left[h\left(z\right)-\textswab{X}^{2}\left(z\right)\right]^{2}+\left[\textswab{X}'\left(z\right)\right]^{2}\not\equiv 0$, as assumed.
\end{proof}

\begin{proof}[Proof of \textup{\fullref{Corollary}{Cor2C}}]
Equations (\ref{qmpC}) and (\ref{PhiSqrtC}) follow immediately from (\ref{qmp}) and (\ref{PhiSqrt}) by the replacement $x\mapsto z$ and $\Phi(x)\mapsto \textswab{X}(z)$, remembering that this also transforms real Ricatti equation (\ref{Ricatti}) to its complex counterpart  (\ref{RicattiC}).
\end{proof}

\begin{proof}[Proof of \textup{\fullref{Theorem}{LemL2C}}]
Let $\pm i\textswab{X}(z)$  be a solution of the complex Ricatti equation (\ref{RicattiC}). Then
\begin{equation}
    \pm i \textswab{X}'\left(z\right)-\textswab{X}^2\left(z\right)+h\left(z\right)=0.
\end{equation}
Since $h(z)\neq \textswab{X}^2(z)$, $\textswab{X}'\left(z\right)$ never vanishes and consequently two sign possibilities are mutually exclusive, i.e. either $i\textswab{X}(z)$ or $-i\textswab{X}(z)$ solves (\ref{RicattiC}).

 On the other hand, under the above assumption, from (\ref{qmpL}) we find that 
\begin{equation}
     \Theta_{\mathtt{top}}\left(z\right)= \Theta_{\mathtt{bot}}\left(z\right)=\pm i\textswab{X}\left(z\right).
\end{equation}
Therefore, $\Theta_{\mathtt{top}}\left(z\right)$ and $\Theta_{\mathtt{bot}}\left(z\right)$ are as well solutions of (\ref{RicattiC}). According to \fullref{Corollary}{Cor2C}, $\textswab{X}\left(z\right)$ is a solution of (\ref{ELexpC}).
\end{proof}

\subsection{Proofs for K{\"a}hler-Norden geometry}\label{proofs2.3}

\begin{proof}[Proof of \textup{\fullref{Lemma}{Einstein}}]
As follows from \fullref{Lemma}{lem-hypC}, $G_h$ is Einstein with  $\eta=-1$. Due to Theorem 5.1 from \cite{Borowiec} (see \cite{Olszak} for a generalization), $\hat{\mathbf{g}}_{h}$ would be Einstein with the same constant if defined according to Eq. (\ref{correpon}) with $\alpha=2$. Since we have chosen $\alpha=1$, by rescaling we get $\eta=-2$. The value of the Ricci scalar simply follows, bearing in mind that we work with a four-dimensional manifold.
\end{proof}

\begin{proof}[Proof of \textup{\fullref{Theorem}{GeoInherited}}]
The convexity requirement put on the set $\mathcal{V}_0$ can always be imposed locally, in a neighborhood of the point $\upsilon=0$. It serves the purpose to assure that the path $\gamma_z^\circ$ given by $\upsilon(s)= s$ belongs to   $\mathcal{V}_0$. If this holds, then (\ref{integralGamma}) reduces to
\begin{equation}\label{integralGammaproof}
    \int_{\gamma_z}\!\!d\zeta\, \Theta_{\mathtt{top}/\mathtt{bot}}(\zeta)=\int_0^1 \!\!ds\,\frac{d\zeta\left(s\right)}{d s}\Theta_{\mathtt{top}/\mathtt{bot}}\left(\zeta\left(s\right)\right),
\end{equation}
being a functional of $\textswab{X}\left(z(s)\right)$. The latter is the geodesic curve in explicit form related with a geodesic curve in $\mathbb{M}_h^\mathbb{C}$, for which the holomorphic parameter $\upsilon$ is restricted to the
real axis. Such a holomorphic curve is also a smooth (real analytic)  geodesic curve on the underlying smooth pseudo-Riemannian manifold \cite{LeBrunPhD,LeBrun2}. Locally, such a curve $\mathbb{R}\ni s \mapsto \left(x(s),\varPhi(s),y(s),\varPsi(s)\right)\in \mathcal{H}_\mathrm{KN}$ obeys the geodesic equation with the Christoffel symbols $\hat{\Gamma}_{\:\:\:jk}^{i}$, inherited from $\Upsilon_{\:\:\:jk}^{i}$ explicitly provided in (\ref{ChristoCom}), which are as follows (see §1.1 of \cite{LeBrunPhD}):
\begin{subequations}\label{ChristoKN}
\begin{equation}
\mathrm{Re}\left[\Upsilon_{\:\:\:zz}^{z}\right]=\hat{\Gamma}_{\:\:\:xx}^{x}=\hat{\Gamma}_{\:\:\:xy}^{y}=-\hat{\Gamma}_{\:\:\:yy}^{x},
\end{equation}
\begin{equation}
\mathrm{Im}\left[\Upsilon_{\:\:\:zz}^{z}\right]=\hat{\Gamma}_{\:\:\:xx}^{y}=-\hat{\Gamma}_{\:\:\:xy}^{x}=-\hat{\Gamma}_{\:\:\:yy}^{y},
\end{equation}
\begin{equation}
\mathrm{Re}\left[\Upsilon_{\:\:\:z\mathfrak{X}}^{z}\right]=\hat{\Gamma}_{\:\:\:x\varPhi}^{x}=\hat{\Gamma}_{\:\:\:x\varPsi}^{y}=-\hat{\Gamma}_{\:\:\:y\varPsi}^{x},
\end{equation}
\begin{equation}
\mathrm{Im}\left[\Upsilon_{\:\:\:z\mathfrak{X}}^{z}\right]=\hat{\Gamma}_{\:\:\:x\varPhi}^{y}=-\hat{\Gamma}_{\:\:\:x\varPsi}^{x}=-\hat{\Gamma}_{\:\:\:y\varPsi}^{y},
\end{equation}
\begin{equation}
\mathrm{Re}\left[\Upsilon_{\:\:\:zz}^{\mathfrak{X}}\right]=\hat{\Gamma}_{\:\:\:xx}^{\varPhi}=\hat{\Gamma}_{\:\:\:xy}^{\varPsi}=-\hat{\Gamma}_{\:\:\:yy}^{\varPhi},
\end{equation}
\begin{equation}
\mathrm{Im}\left[\Upsilon_{\:\:\:zz}^{\mathfrak{X}}\right]=\hat{\Gamma}_{\:\:\:xx}^{\varPsi}=-\hat{\Gamma}_{\:\:\:xy}^{\varPhi}=-\hat{\Gamma}_{\:\:\:yy}^{\varPsi},
\end{equation}
\begin{equation}
\mathrm{Re}\left[\Upsilon_{\:\:\:\mathfrak{X}\mathfrak{X}}^{\mathfrak{X}}\right]=\hat{\Gamma}_{\:\:\:\varPhi\varPhi}^{\varPhi}=\hat{\Gamma}_{\:\:\:\varPhi\varPsi}^{\varPsi}=-\hat{\Gamma}_{\:\:\:\varPsi\varPsi}^{\varPhi},
\end{equation}
\begin{equation}
\mathrm{Im}\left[\Upsilon_{\:\:\:\mathfrak{X}\mathfrak{X}}^{\mathfrak{X}}\right]=\hat{\Gamma}_{\:\:\:\varPhi\varPhi}^{\varPsi}=-\hat{\Gamma}_{\:\:\:\varPhi\varPsi}^{\varPhi}=-\hat{\Gamma}_{\:\:\:\varPsi\varPsi}^{\varPsi}.
\end{equation}
\end{subequations}
I omit vanishing coefficients and symmetries of the Christoffel symbols.

While guaranteed by construction \cite{LeBrun2}, one can confirm by direct inspection that (\ref{ChristoKN}) are Christoffel symbols of $\hat{\mathbf{g}}_{h}$. Therefore, the geodesic equation on $\mathbb{M}_h^\mathrm{KN}$ simplifies to the form
\begin{eqnarray}
-\ddot{\varPhi} & = & \hat{\Gamma}_{\:\:\:xx}^{\varPhi}\dot{x}^{2}+\hat{\Gamma}_{\:\:\:yy}^{\varPhi}\dot{y}^{2}+2\hat{\Gamma}_{\:\:\:xy}^{\varPhi}\dot{x}\dot{y}+\hat{\Gamma}_{\:\:\:\varPhi\varPhi}^{\varPhi}\dot{\varPhi}^{2}+\hat{\Gamma}_{\:\:\:\varPsi\varPsi}^{\varPhi}\dot{\varPsi}^{2}+2\hat{\Gamma}_{\:\:\:\varPhi\varPsi}^{\varPhi}\dot{\varPhi}\dot{\varPsi}\nonumber\\
 & = & \textrm{Re}\left[\Upsilon_{\:\:\:zz}^{\mathfrak{X}}\dot{z}^{2}+\Upsilon_{\:\:\:\mathfrak{X}\mathfrak{X}}^{\mathfrak{X}} \dot{\mathfrak{X}}^2\right],
\end{eqnarray}
\begin{eqnarray}
-\ddot{\varPsi} & = & \hat{\Gamma}_{\:\:\:xx}^{\varPsi}\dot{x}^{2}+\hat{\Gamma}_{\:\:\:yy}^{\varPsi}\dot{y}^{2}+2\hat{\Gamma}_{\:\:\:xy}^{\varPsi}\dot{x}\dot{y}+\hat{\Gamma}_{\:\:\:\varPhi\varPhi}^{\varPsi}\dot{\varPhi}^{2}+\hat{\Gamma}_{\:\:\:\varPsi\varPsi}^{\varPsi}\dot{\varPsi}^{2}+2\hat{\Gamma}_{\:\:\:\varPhi\varPsi}^{\varPsi}\dot{\varPhi}\dot{\varPsi}\nonumber\\
 & = & \textrm{Im}\left[\Upsilon_{\:\:\:zz}^{\mathfrak{X}}\dot{z}^{2}+\Upsilon_{\:\:\:\mathfrak{X}\mathfrak{X}}^{\mathfrak{X}} \dot{\mathfrak{X}}^2\right],
\end{eqnarray}
and similarly for $\ddot{x}$ and $\ddot{y}$. Consequently
\begin{equation}
   \ddot{\mathfrak{X}} =\ddot{\varPhi}+i \ddot{\varPsi}=-\Upsilon_{\:\:\:zz}^{\mathfrak{X}}\dot{z}^{2}-\Upsilon_{\:\:\:\mathfrak{X}\mathfrak{X}}^{\mathfrak{X}} \dot{\mathfrak{X}}^2,
\end{equation}
reproduces the geodesic equation on $\mathbb{M}_h^\mathbb{C}$.

In the final step we recognize that since $\textswab{X}(z)$ is holomorphic,  $\varPhi(s)=\textrm{Re}\left[\textswab{X}\left(z(s)\right)\right]$ and $\varPsi(s)=\textrm{Im}\left[\textswab{X}\left(z(s)\right)\right]$
must both be functions of $x(s)$ and $y(s)$ which satisfy the Cauchy-Riemann equations in these variables. Conversely, a direct substitution of $\varPhi(s)= \Phi\left(x(s),y(s)\right)$ and $\varPsi(s)= \Psi\left(x(s),y(s)\right)$ into the geodesic equation on $\mathbb{M}_h^\mathrm{KN}$, together with the Cauchy-Riemann equations, leads to Eq. (\ref{ELexpC}) under the assumption $\dot{x}(s)+i \dot{y}(s)\neq0$.
\end{proof}

\subsection{Proofs repeated from \cite{Rudnicki}}\label{Sec5}
For the reader's convenience we repeat a few proofs from \cite{Rudnicki}, those 
relevant for \fullref{Sec.}{subsechyp}. However, we do not prove \fullref{Lemma}{lem-hyp}, since we prove \fullref{Lemma}{LemL1} with all details. Note that we often drop the arguments of the functions, therefore, it is important to remember that all discussed theorems deal with the real independent variable $x$.

\begin{proof}[Proof of \textup{\fullref{Proposition}{Prop1}}]
We start from the substitution $\varPhi\left(s\right)=\Phi\left(x\left(s\right)\right)$,
which leads to:\begin{subequations}
\begin{equation}\label{primeph}
\dot{\varPhi}\left(s\right)=\Phi'\left(x\left(s\right)\right)\dot{x}\left(s\right),
\end{equation}
\begin{equation}
\ddot{\varPhi}\left(s\right)=\Phi'\left(x\left(s\right)\right)\ddot{x}\left(s\right)+\Phi''\left(x\left(s\right)\right)\left[\dot{x}\left(s\right)\right]^{2}.
\end{equation}\end{subequations}
From now on we drop the arguments $s$ and $x(s)$, recalling that $\Phi$ is a function of $x$ while $\varPhi$ is a function of $s$, as well as that the prime denotes the derivative with respect to $x$ while the dot is the derivative with respect to $s$. We substitute $\ddot x$ and $\ddot \varPhi$ from the geodesic equations  getting
\begin{equation}
\frac{\Phi^{4}-h^{2}}{\Phi}\dot{x}^{2}+\frac{\dot{\varPhi}^{2}}{\Phi}=\Phi'\left(\frac{h'}{\Phi^{2}-h}\dot{x}^{2}-2\frac{\Phi^{2}-h}{\Phi\left(\Phi^{2}-h\right)}\dot{x}\dot{\varPhi}\right)+\Phi''\dot{x}^{2}.
\end{equation}
After substituting $\dot{\varPhi}$ from (\ref{primeph}) we collect together two terms with $\left[\Phi'\right]^{2}$. After moving all the terms to the right hand side we get
\begin{equation}\label{finalhelpf}
0=\dot{x}^2\left\{\Phi''-\frac{\Phi^{4}-h^{2}}{\Phi}+\frac{h'\Phi'}{\left(\Phi^{2}-h\right)}-\frac{\left(3\Phi^{2}+h\right)\left[\Phi'\right]^{2}}{\Phi\left(\Phi^{2}-h\right)}\right\}.
\end{equation}
As long as $\dot{x}\neq0$, this is the same equation as Eq. (\ref{ELexp}). By direct inspection we can confirm that, assuming  $\dot{x}\neq0$ and using the Christoffel symbols (\ref{ChristoOrig}), Eq. (\ref{finalhelpf}) can be rewritten to the form
\begin{equation}
\Phi''=-\Gamma_{\:\:\:xx}^{\varPhi}+\Gamma_{\:\:\:xx}^{x}\Phi'+\left(2\Gamma_{\:\:\:\varPhi x}^{x}-\Gamma_{\:\:\:\varPhi\varPhi}^{\varPhi}\right)\left[\Phi'\right]^{2}.
\end{equation}
\end{proof}

\begin{proof}[Proof of \textup{\fullref{Theorem}{thm-main2}}]
The same way as in the former proof  we omit the dependence on the variable $x$. The same applies to proof of  \fullref{Corollary}{Cor2}. In order to simplify the derivation we introduce a norm of a velocity of a geodesic curve in explicit form \begin{subequations}
\begin{equation}
\mathcal{L}=\sqrt{\left(h-\Phi^{2}\right)^{2}+\left[\Phi'\right]^{2}},
\end{equation}
and denote
\begin{equation}\label{Thetapm}
\Theta_{\pm}=\Phi\frac{\Phi'\pm\mathcal{L}}{h-\Phi^{2}}.
\end{equation}
We find
\begin{equation}
\Theta_{\pm}'=\Phi'\frac{\Phi'\pm\mathcal{L}}{h-\Phi^{2}}+\Phi\frac{\left(\Phi''\pm\mathcal{L}'\right)\left(h-\Phi^{2}\right)-\left(\Phi'\pm\mathcal{L}\right)\left(h'-2\Phi\Phi'\right)}{\left(h-\Phi^{2}\right)^{2}},
\end{equation}
and
\begin{equation}
\Theta_{\pm}^{2}=2\Phi^{2}\Phi'\frac{\Phi'\pm\mathcal{L}}{\left(h-\Phi^{2}\right)^{2}}+\Phi^{2}.
\end{equation}    
\end{subequations}\begin{subequations}
We then get
\begin{equation}\label{RicPQ}
\Theta_{\pm}'+\Theta_{\pm}^{2}+h=P\pm Q,
\end{equation}
where
\begin{equation}
P=\frac{\left[\Phi'\right]^{2}+\Phi\Phi''}{h-\Phi^{2}}-\Phi\Phi'\frac{h'-4\Phi\Phi'}{\left(h-\Phi^{2}\right)^{2}}+\Phi^{2}+h,
\end{equation}
\begin{equation}
Q=\frac{\Phi'\mathcal{L}+\Phi\mathcal{L}'}{h-\Phi^{2}}-\Phi\frac{h'-4\Phi\Phi'}{\left(h-\Phi^{2}\right)^{2}}\mathcal{L}.
\end{equation}\end{subequations}

We first rearrange\begin{subequations}
\begin{eqnarray}
P & = & \frac{\left[\Phi'\right]^{2}+\Phi\Phi''}{h-\Phi^{2}}-\Phi\Phi'\frac{h'-4\Phi\Phi'}{\left(h-\Phi^{2}\right)^{2}}+\Phi^{2}+h\nonumber \\
 & = & \frac{\Phi\Phi''}{h-\Phi^{2}}+\Phi'\frac{\Phi'\left(h-\Phi^{2}\right)-\Phi\left(h'-4\Phi\Phi'\right)}{\left(h-\Phi^{2}\right)^{2}}+\Phi^{2}+h\nonumber \\
 & = & \frac{\Phi\Phi''}{h-\Phi^{2}}+\Phi'\frac{\Phi'\left(h+3\Phi^{2}\right)-\Phi h'}{\left(h-\Phi^{2}\right)^{2}}+\Phi^{2}+h,\label{Pfinal}
\end{eqnarray}
and
\begin{equation}
Q=\frac{\Phi}{h-\Phi^{2}}\mathcal{L}'+\frac{\Phi'\left(h+3\Phi^{2}\right)-\Phi h'}{\left(h-\Phi^{2}\right)^{2}}\mathcal{L}.
\end{equation}\end{subequations}
Next, since\begin{subequations}
\begin{equation}
\mathcal{L}'=\frac{\left(h-\Phi^{2}\right)\left(h'-2\Phi\Phi'\right)+\Phi'\Phi''}{\mathcal{L}},
\end{equation}
we observe that $Q/\Phi'$ is regular at $\Phi'= 0$. Therefore, we can write
\begin{eqnarray}
\frac{\mathcal{L}Q}{\Phi'} & = & \frac{\Phi}{h-\Phi^{2}}\frac{\mathcal{L}\mathcal{L}'}{\Phi'}+\frac{h\Phi'+\Phi\left(3\Phi\Phi'-h'\right)}{\left(h-\Phi^{2}\right)^{2}}\frac{\mathcal{L}^{2}}{\Phi'}\nonumber \\
 & = & \frac{\Phi\Phi''}{h-\Phi^{2}}+\frac{h'-2\Phi\Phi'}{\Phi^{-1}\Phi'}+\frac{\Phi'\left(h+3\Phi^{2}\right)-\Phi h'}{\left(h-\Phi^{2}\right)^{2}}\left\{ \frac{\left(h-\Phi^{2}\right)^{2}}{\Phi'}+\Phi'\right\} \nonumber \\
 & = & \frac{\Phi\Phi''}{h-\Phi^{2}}+\frac{\bcancel{h'}-\cancel{2\Phi\Phi'}}{\xcancel{\Phi^{-1}\Phi'}}+\left\{ h+\cancel{3}\Phi^{2}-\xcancel{\frac{\Phi}{\Phi'}}\bcancel{h'}+\frac{\Phi'\left(h+3\Phi^{2}\right)-\Phi h'}{\left(h-\Phi^{2}\right)^{2}}\Phi'\right\} \nonumber \\
 & = & \frac{\Phi\Phi''}{h-\Phi^{2}}+h+\Phi^{2}+\frac{\Phi'\left(h+3\Phi^{2}\right)-\Phi h'}{\left(h-\Phi^{2}\right)^{2}}\Phi'\nonumber \\ \label{QP}
 & = & P,
\end{eqnarray}
where we used
\begin{eqnarray}
\frac{\Phi}{h-\Phi^{2}}\frac{\mathcal{L}\mathcal{L}'}{\Phi'}&=&\frac{\Phi\left[\left(h-\Phi^{2}\right)\left(h'-2\Phi\Phi'\right)+\Phi'\Phi''\right]}{\Phi'\left(h-\Phi^{2}\right)}\\ \nonumber
& =&\frac{\Phi\Phi''}{h-\Phi^{2}}+\frac{h'-2\Phi\Phi'}{\Phi^{-1}\Phi'},
\end{eqnarray}
while passing from first to second line.  The last equality (\ref{QP}) follows from (\ref{Pfinal}). \end{subequations} 
Just to make it clear, Eq. (\ref{QP}) in fact proves $Q=\Phi' P/\mathcal{L}$, regardless of whether a possibility $\Phi'=0$ occurs or not. We have divided by $\Phi'$ only to make the derivation slightly more compact.

By virtue of (\ref{QP}) the equation (\ref{RicPQ}) becomes
\begin{equation}
\Theta_{\pm}'+\Theta_{\pm}^{2}+h=\left(1\pm\frac{\Phi'}{\mathcal{L}}\right)P.
\end{equation}
To finalize the proof we observe that $1\pm\Phi'/\mathcal{L}$ cannot
vanish. Moreover
\begin{equation}\label{finalproof}
u_{\mathtt{top}/\mathtt{bot}}''+h\, u_{\mathtt{top}/\mathtt{bot}}=\left(\Theta_{\mp}'+\Theta_{\mp}^{2}+h\right)u_{\mathtt{top}/\mathtt{bot}}=\left(1\mp\frac{\Phi'}{\mathcal{L}}\right)P\,u_{\mathtt{top}/\mathtt{bot}},
\end{equation}
so that both $u_{\mathtt{top}}$ and $u_{\mathtt{bot}}$ are solutions
of (\ref{Eq.main}) if and only if $P=0$. Consequently, $u = A u_{\mathtt{top}} +B u_{\mathtt{bot}}$, with arbitrary constants $A$ and $B$, is the general solution of (\ref{Eq.main}) if and only if $P=0$. In the final step we make a simple rearrangement
\begin{eqnarray}
P & = & \frac{\Phi\Phi''}{h-\Phi^{2}}+\Phi'\frac{\Phi'\left(h+3\Phi^{2}\right)-\Phi h'}{\left(h-\Phi^{2}\right)^{2}}+\Phi^{2}+h\nonumber \\
 & = & \frac{\Phi}{h-\Phi^{2}}\left[\Phi''-\Phi'\frac{\Phi'\left(h+3\Phi^{2}\right)-\Phi h'}{\Phi\left(\Phi^{2}-h\right)}-\frac{\Phi^{4}-h^{2}}{\Phi}\right].\label{Prearfin}
\end{eqnarray}
to recognize that $P=0$ if and only if (\ref{ELexp}) holds.
\end{proof}

\begin{proof}[Proof of  \textup{\fullref{Corollary}{Cor2}}]
By comparison between (\ref{qmp}) and (\ref{Thetapm}) we recognize that $\Theta_-=\Theta_\mathtt{top}$ and $\Theta_+=\Theta_\mathtt{bot}$. Therefore, the first assertion of \fullref{Corollary}{Cor2} follows directly from (\ref{finalproof}) and \fullref{Theorem}{thm-main2}. The second assertion given as Eq. (\ref{PhiSqrt}) follows immediately from multiplication of (\ref{qmp1}) and (\ref{qmp2}).
\end{proof}

\section{Acknowledgements}
\label{sec:acks}
The author is grateful to Jacek Gulgowski for a helpful discussion, and Claude LeBrun for valuable correspondence.

\bibliography{aomart2}
\bibliographystyle{aomalpha}

\end{document}